\documentclass[a4paper, 12pt]{amsart}
\usepackage{amssymb, hyperref,amsmath}

\newcommand\field[1]{\mathbb{#1}}
\newcommand\NN{\field{N}}

\newcommand\TT{\field{T}}
\newcommand\ZZ{\field{Z}}
\newcommand\Bb{\mathcal B}

\newcommand\Gg{\mathcal G}
\newcommand\Hh{\mathcal H}

\renewcommand\ker{\operatorname{ker}}

\newcommand\Aut{\operatorname{Aut}}

\newcommand\MCE{\operatorname{MCE}}

\newcommand\Ind{\operatorname{Ind}}
\newcommand\supp{\operatorname{supp}}

\theoremstyle{plain}
\newtheorem{theorem}{Theorem}[section]
\newtheorem{cor}[theorem]{Corollary}
\newtheorem{lemma}[theorem]{Lemma}
\newtheorem{prop}[theorem]{Proposition}

\theoremstyle{remark}
\newtheorem{rmk}[theorem]{Remark}

\theoremstyle{definition}
\newtheorem{dfn}[theorem]{Definition}

\newtheorem{notation}[theorem]{Notation}

\numberwithin{equation}{section}

\title{Uniqueness theorems for topological higher-rank graph $C^*$-algebras}
\author{Jean Renault}
\address{D\'epartement de Math\'ematiques\\
Universit\'e d'Orl\'eans\\
BP 6759\\
45067 Orl\'eans Cedex 2\\
FRANCE}
\email{Jean.Renault@univ-orleans.fr}

\author{Aidan Sims}
\address{School of Mathematics and Applied Statistics\\
Austin Keane Building (15)\\
University of Wollongong\\
NSW 2522\\
AUSTRALIA}
\email{asims@uow.edu.au}

\author{Dana P. Williams}
\address{Department of Mathematics \\
Dartmouth College \\
Hanover, NH 03755-3551\\
USA}
\email{dana.williams@Dartmouth.edu}

\author{Trent Yeend}
\address{School of Mathematical and Physical Sciences\\
Building~V\\
University of Newcastle\\
Callaghan NSW 2308\\
AUSTRALIA}
\email{Trent.Yeend@ihpa.gov.au}

\date{\today}
\thanks{This research was supported by the Australian Research Council.}

\begin{document}

\begin{abstract}
We consider the boundary-path groupoids of topological higher-rank graphs. We show that the all
such groupoids are topologically amenable. We deduce that the $C^*$-algebras of topological
higher-rank graphs are nuclear and prove versions of the gauge-invariant uniqueness theorem and the
Cuntz-Krieger uniqueness theorem. We then provide a necessary and sufficient condition for
simplicity of a topological higher-rank graph $C^*$-algebra, and a condition under which it is also
purely infinite.
\end{abstract}

\maketitle

\section{Introduction}

Groupoids are a powerful and widely applicable model for operator algebras. One area of
operator-algebra theory in which they have been particularly prominent recently is the field of
graph algebras and their analogues.

The inception of the field of graph $C^*$-algebras goes back to the work of Cuntz and Krieger
\cite{CK1980}, and the subsequent work of Enomoto and Watatani \cite{EW1980}, on simple purely
infinite $C^*$-algebras associated to finite binary matrices. However, the theory of graph
$C^*$-algebras really took off after the work of Kumjian-Pask-Raeburn-Renault \cite{KPRR1997}. The
analysis there was facilitated by realising the $C^*$-algebras of interest as groupoid
$C^*$-algebras and employing Renault's structure theory~\cite{Renault1980}.

Since then, the class of graph $C^*$-algebras has been generalised in various directions, including
ultragraph $C^*$-algebras \cite{Tomforde2003c}, higher-rank graph $C^*$-algebras \cite{KP2000} and
topological graph $C^*$-algebras \cite{Katsura2004} to name a few. Though these generalisations
have not all been developed using groupoid methods, in each case a natural groupoid model exists
\cite{FMY2005, pp_MM2007, Paterson2002, Yeend2006}.

In 2005, Yeend developed the notion of a topological higher-rank graph, simultaneously
generalising Katsura's notion of a topological graph and Kumjian and Pask's notion of a
higher-rank graph. Yeend associated to each topological higher-rank graph $\Lambda$ a
groupoid $\Gg_\Lambda$ and hence a $C^*$-algebra $C^*(\Lambda) := C^*(\Gg_\Lambda)$.
Yeend's construction is sufficiently general to capture Katsura's algebras and the
finitely aligned $k$-graph $C^*$-algebras of \cite{RSY2004}. However, the question of
amenability of $\Gg_\Lambda$ remained unresolved in general, so Yeend's key
$C^*$-algebraic results held only under additional hypotheses. In addition, the
injectivity hypothesis on Yeend's uniqueness theorems is phrased in terms of functions on
$\Gg_\Lambda$ rather than in terms of the underlying topological $k$-graph $\Lambda$.

The first version of the current paper, posted by Renault, Sims and Yeend on the arXiv preprint
server in 2009, aimed to resolve the question of amenability of Yeend's groupoid and to prove
versions of the gauge-invariant uniqueness theorem and the Cuntz-Krieger uniqueness theorem with an
injectivity hypothesis involving only the algebra of continuous functions on the vertex set of the
topological graph. Our approach to amenability was to show that the kernel of the canonical
$\ZZ^k$-valued cocycle $c$ on $\Gg_\Lambda$ was amenable (by providing a measure-theoretic
direct-limit decomposition into equivalence relations) and then bootstrap up to amenability of
$\Gg_\Lambda$ by composing invariant means on $c^{-1}(0)$ with the mean on $\ZZ^k$. Shortly after
submission, Williams spotted an error in our bootstrapping argument. The paper was withdrawn and
Williams became involved as we considered how to repair the gap. In the mean time, versions of the
key results, namely the gauge-invariant uniqueness theorem and the Cuntz-Krieger uniqueness
theorem, were proved respectively in \cite{CLSV} and \cite{Yamashita:xx09} using the machinery of
product systems (though our Proposition~\ref{prp:kernel amenable} is required to identify the
$C^*$-algebras described there with Yeend's $C^*$-algebra). As a result, this project lay dormant
for some time.

Two recent developments brought the project out of mothballs. The first is Spielberg's clever
argument \cite[Proposition~9.3]{Spielberg:xx11} which combines groupoid theory and coaction theory
to show that if $c$ is a cocycle from an \'etale Hausdorff groupoid $\Gg$ into a countable abelian
group $G$ and $c^{-1}(0)$ is amenable, then $\Gg$ is amenable; this fixed the gap in our original
argument. The second is the recent characterisation of simplicity for the $C^*$-algebra of a
second-countable locally compact Hausdorff \'etale amenable groupoid $\Gg$ in
\cite{BrownClarkEtAl:xx12}: Yeend showed in \cite{Y1} that $\Gg_\Lambda$ has all these properties
except for amenability, and simplicity was not addressed in either \cite{CLSV} or
\cite{Yamashita:xx09}.

In this revised article, we combine Spielberg's argument with our previous analysis to prove that
$\Gg_\Lambda$ is amenable and prove a gauge-invariant uniqueness theorem. The proof of what is now
Proposition~\ref{prp:Lambda^0->G^0} has been significantly simplified by
\cite[Proposition~5.16]{CLSV}. We then use the results of \cite{BrownClarkEtAl:xx12} to prove a
version of the Cuntz-Krieger uniqueness theorem and to provide a necessary and sufficient condition
for simplicity of $C^*(\Lambda)$. We conclude by providing a sufficient condition for
$C^*(\Lambda)$ to be purely infinite.

\subsection*{Acknowledgement}
We thank Toke Carlsen for his reading of a preliminary draft of the manuscript, and for his
valuable comments and input.

\section{Background}\label{sec:prelims}
Our results require a fair amount of background. We do not give full detail, especially
as regards the theory of groupoids. For more detail see, for example, \cite{A-DR,
Renault1980, Y1}.

We regard $\NN^k$ as a semigroup with identity $0$, or sometimes as a category with a single object
and composition defined by the addition operation. For $m,n \in \NN^k$, we say $m \le n$ if $m_i
\le n_i$ for all $i \in \{1, \dots, k\}$. We write $m \vee n$ for the coordinatewise maximum of $m$
and $n$. We frequently work also with the set $(\NN \cup \{\infty\})^k$; we extend the addition
operation and the order $\le$ from $\NN^k$ to $(\NN \cup \{\infty\})^k$ in the obvious way.

\subsection{Groupoids}

For details of what follows, see \cite{Renault1980}. A groupoid $\Gg$ is a small category with
inverses. We denote the domain and codomain maps by $s$ and $r$, the unit space by $\Gg^{(0)}$, and
the collection of composable pairs by $\Gg^{(2)}$. A topological groupoid is a groupoid endowed
with a topology under which both the inversion map and the composition map are continuous. In this
paper, we consider only locally compact Hausdorff groupoids. A groupoid is \'etale if the range map
is a local homeomorphism; when the topology is Hausdorff, it follows that $\Gg^{(0)}$ is both open
and closed in $\Gg$. For each unit $u$ of an \'etale groupoid $\Gg$ the sets $\Gg^u := r^{-1}(u)$
and $\Gg_u := s^{-1}(u)$ are discrete.

The \emph{isotropy} at a unit $u \in \Gg^{(0)}$ is the group $\Gg^u_u : r^{-1}(u) \cap s^{-1}(u)$.
The phrase ``points with trivial isotropy" is frequently used in the literature to refer to the
units $u$ of a groupoid $\Gg$ such that the isotropy at $u$ is the trivial group. We say that a
groupoid $\Gg$ is \emph{principal} if every unit $u$ has trivial isotropy; algebraically, these are
just equivalence relations, but topologically they can be quite different.

The \emph{orbit} of a unit $u$ of a groupoid $\Gg$ is the set
\[
[u] = \{r(x) : s(x) = u\};
\]
that is, $[u] = r(\Gg_u) = s(\Gg^u)$. A subset $U$ of $\Gg^{(0)}$ is \emph{invariant} if $[u]
\subseteq U$ whenever $u \in U$.

\subsection{Groupoid {\texorpdfstring{{$C^*$}}{{C*}}}-algebras}

We will now summarise the constructions of the full- and reduced groupoid $C^*$-algebras of an
\'etale locally compact Hausdorff groupoid $\Gg$. These constructions can be carried through for
any groupoid admitting a Haar system, but the formulae are simpler in our situation. For full
details see \cite[Section~II.1]{Renault1980}; or for a detailed treatment of \'etale groupoids, see
\cite{Paterson:Groupoidsinversesemigroups99} or \cite[Section~3]{Exel2008}.

Consider the space $C_c(\Gg)$ of compactly supported complex-valued functions on $\Gg$. For $x \in
\Gg$ and $f \in C_c(\Gg)$, the set $\{y : r(y) = r(x), f(y) \not= 0\}$ is both compact and discrete
and hence finite. Thus we may sensibly define an operation $\ast : C_c(\Gg) \times C_c(\Gg) \to
C_c(\Gg)$ by
\[
(f \ast g)(x) = \sum_{r(y) = r(x)} f(y) g(y^{-1}x).
\]
The space $C_c(\Gg)$ becomes a topological $^*$-algebra with the involution $f^*(x) =
\overline{f(x^{-1})}$ and the convolution product $\ast$ defined above.

A representation of $C_c(\Gg)$ is a nondegenerate $^*$-homomorphism $\pi : C_c(\Gg) \to B(\Hh)$
which is continuous from the inductive limit topology on $C_c(\Gg)$ to the strong-operator topology
on $\Bb(\Hh)$. Renault's disintegration theorem \cite[Proposition~4.2]{Renault1987} together with
\cite[Propositions II.1.7~and~II.1.11]{Renault1980} implies that there is a pre-$C^*$-norm on
$C_c(\Gg)$ determined by
\[
\|f\| = \sup\{ \|\pi(f)\| : \pi\text{ is a representation of $C_c(\Gg)$}\}.
\]
The full groupoid $C^*$-algebra $C^*(\Gg)$ is the $C^*$-completion of $C_c(\Gg)$ with respect to
this norm.

To define the reduced groupoid $C^*$-algebra, fix $u \in \Gg^{(0)}$ and let $\ell^2(\Gg_u)$ be the
Hilbert space with orthonormal basis $\{\xi_x : x \in \Gg_u\}$. There is a representation $\Ind
\epsilon_u : C_c(\Gg) \to \Bb(\ell^2(\Gg_u))$ such that for $f \in C_c(\Gg)$ and $x \in \Gg_u$, we
have
\[
\Ind \epsilon_u(f) \xi_x = \sum_{y \in \Gg_u} f(yx^{-1}) \xi_y.
\]
The reduced groupoid $C^*$-algebra $C^*_r(\Gg)$ is then the completion of $C_c(\Gg)$ in the
$C^*$-norm
\[
\|f\|_r = \sup_{u \in \Gg^{(0)}} \|\Ind \epsilon_u(f)\|.
\]

There are at least two notions of amenability for groupoids: (topological) amenability
\cite[Definition~2.2.8]{A-DR}, and measurewise amenability \cite[Definition~3.3.1]{A-DR}.
Topological amenability of $\Gg$ implies measurewise amenability of $\Gg$ \cite[page~83]{A-DR}
which in turn implies that $C^*(\Gg)$ and $C^*_r(\Gg)$ coincide \cite[Proposition~6.18]{A-DR}. If
$\Gg$ is both \'etale and second-countable, then it has countable orbits and a continuous Haar
system (consisting of counting measures), and in this case \cite[Theorem~3.3.7]{A-DR} implies that
measurewise and topological amenability are equivalent.

\subsection{Topological higher-rank graphs}

For the details of this and the next section, see \cite{YeendPhD, Y1}. A $k$-graph is a small
category $\Lambda$ equipped with a functor $d : \Lambda \to \NN^k$ which satisfies the
factorisation property: if $d(\lambda) = m+n$ then there exist unique $\mu \in d^{-1}(m)$ and $\nu
\in d^{-1}(n)$ such that $\lambda = \mu\nu$. We call $d$ the \emph{degree map} on $\Lambda$, and
denote $d^{-1}(n)$ by $\Lambda^n$.\footnote{When $k = 1$ so that $n \in \NN$, there is a slight
clash of notation here with the usual notation for the product space $\prod^n_{i=1} \Lambda$; but
the meaning is usually clear from context.} An argument using the factorisation property shows that
$\Lambda^0$ is equal to the collection of identity morphisms of $\Lambda$. So the domain and
codomain functions determine maps $s,r : \Lambda \to \Lambda^0$ which we call the source and range
maps.

We regard a $k$-graph $\Lambda$ as a kind of generalised directed graph: we think of $\Lambda^0$ as
a collection of vertices; we think of each $\lambda \in \Lambda$ as a path from $s(\lambda)$ to
$r(\lambda)$; and the degree map $d$ plays the role of a generalised length function.

Given sets $X,Y \subseteq \Lambda$ we write $XY$ for the set $\{\mu\nu : \mu \in X, \nu
\in Y, s(\mu) = r(\nu)\}$. In particular, for $V \subseteq \Lambda^0$ and $X \subseteq
\Lambda$, $VX = \{\lambda \in X : r(\lambda) \in V\}$ and $XV = \{\lambda \in X :
s(\lambda) \in V\}$. By the usual abuse of notation, for a singleton set $\{\lambda\}
\subseteq \Lambda$, we write $\Lambda \lambda$ and $\lambda \Lambda$ in place of
$\Lambda\{\lambda\}$ and $\{\lambda\}\Lambda$.

The factorisation property implies that if $\lambda \in \Lambda$ and $m \le n \le d(\lambda)$, then
there are unique $\lambda(0,m) \in \Lambda^m$, $\lambda(m,n) \in \Lambda^{n-m}$ and
$\lambda(n,d(\lambda)) \in \Lambda^{d(\lambda) - n}$ such that $\lambda =
\lambda(0,m)\lambda(m,n)\lambda(n,d(\lambda))$. We think of $\lambda(m,n)$ as the segment of
$\lambda$ from position $m$ to position $n$ along $\lambda$.

Given $\mu$ and $\nu$ in $\Lambda$, we say that $\lambda$ is a common extension of $\mu$ and $\nu$
if we can factorise $\lambda = \mu\mu'$ and $\lambda = \nu\nu'$ for some $\mu', \nu' \in \Lambda$.
We say that $\lambda$ is a minimal common extension of $\mu$ and $\nu$ if it is a common extension
such that $d(\lambda) = d(\mu) \vee d(\nu)$. We denote by $\MCE(\mu,\nu)$ the set of all minimal
common extensions of $\mu$ and $\nu$. If $r(\mu) \not= r(\nu)$, then $\MCE(\mu,\nu) = \emptyset$.
Given subsets $X,Y \subseteq \Lambda$, we define
\[
\MCE(X,Y) := \bigcup_{\mu \in X, \nu \in Y} \MCE(\mu,\nu).
\]

A topological $k$-graph is a $k$-graph $\Lambda$ endowed with a second-countable locally compact
Hausdorff topology such that each $\Lambda^n$ is open, the range map is continuous, composition is
continuous and open, and the source map is a local homeomorphism.

We say that the topological $k$-graph $\Lambda$ is compactly aligned if, for every pair of compact
subsets $X,Y \subseteq \Lambda$, the set $\MCE(X,Y)$ is also compact. Given $v \in \Lambda^0$ we
say that a subset $E$ of $\Lambda$ is compact exhaustive for $v$ if $E$ is compact, $r(E)$ is a
neighbourhood of $v$, and for all $\lambda \in r(E)\Lambda$ there exists $\mu \in E$ such that
$\MCE(\lambda,\mu) \not= \emptyset$.

An important class of examples of higher-rank graphs, which we will use to make sense of the notion
of an infinite path in a $k$-graph, are the (discrete) higher-rank graphs $\Omega_{k,m}$. These are
defined as follows. For fixed $k \ge 1$ and $m \in (\NN \cup \{\infty\})^k$, the $k$-graph
$\Omega_{k,m}$ has morphisms $\{(p,q) : p,q \in \NN^k, p \le q \le m\}$. The range and source maps
are $r(p,q) = (p,p)$ and $s(p,q) = (q,q)$, composition is determined by $(p,q)(q,r) = (p,r)$, and
the degree map is given by $d(p,q) = q - p$. We usually abbreviate a vertex $(p,p)$ of
$\Omega_{k,m}$ as $p$.

A $k$-graph morphism $x$ from a $k$-graph $\Lambda$ to a $k$-graph $\Gamma$ is a functor $x :
\Lambda \to \Gamma$ which intertwines the degree maps. Given a $k$-graph $\Lambda$, each $\lambda
\in \Lambda$ determines, and is determined by, the unique $k$-graph morphism $x_\lambda :
\Omega_{k, d(\lambda)} \to \Lambda$ such that $x_\lambda(m,n) = \lambda(m,n)$ for all $m \le n \le
d(\lambda)$. By analogy, for arbitrary $m \in (\NN \cup \{\infty\})^k$, we regard a $k$-graph
morphism $x : \Omega_{k,m} \to \Lambda$ as a (possibly infinite) path in $\Lambda$, and define
$r(x) := x(0)$, and $d(x) = m$. If $m_i < \infty$ for all $i$, then we also write $s(x)$ for
$x(m)$, but if $m_i = \infty$ for some $i$, then $x$ has no source.

For $m \in (\NN \cup \{\infty\})^k$, a boundary path of degree $m$ in a topological $k$-graph
$\Lambda$ is a $k$-graph morphism $x : \Omega_{k,m} \to \Lambda$ such that for every $n \in \NN^k$
with $n \le m$ and each compact exhaustive set $E$ for $x(n,n)$, there is an element $\lambda$ of
$E$ such that $d(x) \ge n+d(\lambda)$ and $x(n,n+d(\lambda)) = \lambda$. Fix a boundary path $x$ of
degree $m$. For each $n \in \NN^k$ with $n \le m$, there is a unique boundary path $\sigma^n(x)$ of
degree $m - n$ defined by $\sigma^n(x)(p,q) = x(n+p, n+q)$ for all $(p,q) \in \Omega_{k, m - n}$.
Given $\mu \in \Lambda$ with $s(\mu) = r(x)$, there is a unique boundary path $\mu x$ of degree
$d(\mu) + m$ such that $(\mu x)(0, d(\mu)) = \mu$ and such that $(\mu x)(p + d(\mu), q + d(\mu)) =
x(p,q)$ for all $(p,q) \in \Omega_{k,m}$. We have $\sigma^{d(\mu)}(\mu x) = x = x(0,
n)\sigma^n(x)$.

We denote the collection of all boundary paths in $\Lambda$ by $\partial\Lambda$. For a subset $U$
of $\Lambda$, we denote by $Z(U)$ the collection
\[
\{x \in \partial\Lambda : x(0,n) \in U\text{ for some } n\in\NN^k\text{ with }n \le d(x)\}.
\]
The collection of sets
\begin{equation}\label{eq:topology base}
\{Z(U) \cap Z(F)^c : U \subseteq \Lambda^m\text{ is relatively compact and open, } F \subseteq \overline{U}\Lambda\text{ is compact}\}
\end{equation}
form a basis for a locally compact Hausdorff topology on $\partial\Lambda$.

\subsection{The boundary-path groupoid of a topological
higher-rank graph}

Given a compactly aligned topological $k$-graph $\Lambda$, we define a set $\Gg_\Lambda$ by
\[
\Gg_\Lambda := \{(x, m - n, y) : m,n \in \NN^k, x,y \in \partial\Lambda, m \le d(x), n \le d(y)\text{ and } \sigma^m(x) = \sigma^n(y)\}.
\]
Define $\Gg_\Lambda^{(0)} := \{(x,0,x) : x \in
\partial\Lambda\}$, and identify it with $\partial\Lambda$ via $(x,0,x) \mapsto x$. For $(x,p,y) \in \Gg_\Lambda$, define
$r(x,p,y) = x$ and $s(x,p,y) = y$. With structure maps
\[
(x,p,y)^{-1} = (y,-p,x),\qquad\text{and}\qquad
    (x,p,y)(y,q,z) = (x,p+q,z),
\]
the set $\Gg_\Lambda$ becomes a groupoid with unit space $\partial\Lambda$, and $c(x,p,y) := p$
defines a continuous 1-cocycle $c : \Gg_\Lambda \to \ZZ^k$.

For $U,V \subseteq \Lambda$, define $U \ast_s V := \{(\mu,\nu) \in U \times V : s(\mu)=s(\nu)\}$.
For $F \subseteq \Lambda \ast_s \Lambda$ and $p \in \ZZ^k$, define
\[
Z(F, p) := \{(\mu x, p, \nu x) : (\mu,\nu) \in F, d(\mu) - d(\nu) = p, s(\mu)=s(\nu), x \in s(\mu)\partial\Lambda\}.
\]

The following follows from Yeend's results, but it is worthwhile to state it explicitly.

\begin{lemma}\label{lem:GLambda}
Let $\Lambda$ be a compactly aligned topological $k$-graph. Then
\begin{align*}
\{ Z(U \ast_s V, p-q) \cap Z(F, p-q)^c : {}&p,q \in \NN^k,\  U \subseteq \Lambda^p\text{ and }V \subseteq \Lambda^q, \\
    &\text{$U,V$ are relatively compact and open,}\\
    &\textstyle\text{and }F\text{ is a compact subset of } \bigcup_{\alpha \in \Lambda} \overline{U}\alpha \times \overline{V}\alpha\}
\end{align*}
is a basis for a locally compact Hausdorff topology on $\Gg_\Lambda$ under which $\Gg_\Lambda$
becomes a locally compact Hausdorff \'etale groupoid.
\end{lemma}
\begin{proof}
Proposition~3.6 and Theorem~3.16 of \cite{Y1} imply that the sets of the form described
in the lemma are a basis for a second-countable, locally compact, Hausdorff topology on
the path groupoid $G_\Lambda$ of $\Lambda$, and that $G_\Lambda$ is an \'etale groupoid
under this topology. Propositions 4.4~and~4.7 of \cite{Y1} show that $\partial\Lambda$ is
a closed invariant subset of $G_\Lambda^{(0)}$. Since $\Gg_\Lambda$ is by definition the
restriction of $G_\Lambda$ to $\partial\Lambda$, the result follows.
\end{proof}

\begin{notation}\label{ntn:lazy}
Given $U \subseteq \Lambda^p$ and $V \subseteq\Lambda^q$ and a compact subset $F \subseteq
\bigcup_{\alpha \in \Lambda} \overline{U}\alpha \times \overline{V}\alpha$, it is unambiguous to
abbreviate the basic open set $Z(U *_s V, p-q) \cap Z(F, p-q)^c$ as $Z(U *_s V) \cap Z(F)^c$, and
we will frequently do so.
\end{notation}

The topological higher-rank graph $C^*$-algebra $C^*(\Lambda)$ is defined to be the full groupoid
$C^*$-algebra $C^*(\Gg_\Lambda)$.

\section{Injectivity of representations on functions on the unit space}\label{sec:C0(G0)}

The uniqueness theorems in \cite{Y1} start with a representation of $C^*(\Gg_\Lambda)$ which
restricts to an injection of $C_0(\Gg_\Lambda^{(0)}) = C_0(\partial\Lambda)$. For graph
$C^*$-algebras, topological graph $C^*$-algebras and higher-rank graph $C^*$-algebras, the usual
hypothesis is that the given representation be injective on the embedded copy of $C_0(\Lambda^0)$.
We show that the two hypotheses are equivalent by showing that injectivity on $C_0(\Lambda^0)$
implies injectivity on  $C_0(\partial\Lambda)$. That is, the usual hypothesis also suffices for
topological higher-rank graphs.

The definition of the topology on $\partial\Lambda$ ensures that the range map $r : x \mapsto x(0)$
is continuous from $\partial\Lambda$ to $\Lambda^0$. Proposition~4.3 of \cite{Y1} implies that $r$
is surjective. It therefore induces an injection
\[
r^* : C_0(\Lambda^0) \hookrightarrow C_0(\partial\Lambda)\quad\text{ such that } r^*(f) = f\circ r
\quad\text{ for all $f \in C_0(\Lambda^0)$.}
\]

\begin{prop}\label{prp:Lambda^0->G^0}
Let $\Lambda$ be a compactly aligned topological $k$-graph. Let $\pi$ be a representation
of $C^*(\Gg_\Lambda)$. If $\pi|_{r^*(C_c(\Lambda^0))}$ is injective, then
$\pi|_{C_0(\Gg_\Lambda^{(0)})}$ is injective.
\end{prop}
\begin{proof}
The ideal $\ker(\pi) \cap C_0(\partial\Lambda)$ consists of all functions supported on some open
invariant subset $U$ of $\partial\Lambda$. So $X := \partial\Lambda \setminus U$ is a closed
invariant set and $\pi$ factors through a representation of $\Gg_\Lambda|_X$. That $\pi \circ r$ is
injective on $C_c(\Lambda^0)$ implies that $X \cap Z(V) \not= \emptyset$ for every open $V
\subseteq \Lambda^0$. Fix $v \in \Lambda^0$ and a fundamental sequence of compact neighbourhoods
$(K_n)^\infty_{n=1}$ of $v \in \Lambda^0$. We have just seen that $X \cap K_n \not= \emptyset$ for
all $n$, so fix a sequence $(x_n)^\infty_{n=1}$ with each $x_n \in X \cap K_n$. By compactness we
may pass to a convergent subsequence with limit $x$, and then continuity of the range map ensures
that $r(x) = v$. Hence $K \cap v\partial\Lambda \not= \emptyset$.

Proposition~5.16 of \cite{CLSV} implies that the only closed invariant set of $\partial\Lambda$
which intersects each $v\partial\Lambda$ is $\partial\Lambda$ itself. So $K = \partial\Lambda$, and
hence $\ker(\pi) \cap C_0(\partial\Lambda) = \{0\}$.
\end{proof}

\section{Amenability and the gauge-invariant uniqueness theorem}\label{sec:amenable}

In this section we prove a variant of an Huef and Raeburn's gauge-invariant uniqueness
theorem \cite{HR1997} for topological higher-rank graph $C^*$-algebras. A key ingredient
is amenability of $\Gg_\Lambda$ which guarantees that $C^*(\Gg_\Lambda)$ and
$C^*_r(\Gg_\Lambda)$ coincide; it follows from \cite[Proposition~4.8]{Renault1980} that
the conditional expectation of $C^*(\Gg_\Lambda)$ onto $C_0(\Gg_\Lambda^{(0)})$ is
faithful.

Recall that given a topological higher-rank graph $\Lambda$, we denote by $c$ the canonical
1-cocycle $c : \Gg_\Lambda \to \ZZ^k$ given by $c(x,m,y) = m$.

\begin{theorem}[The gauge-invariant uniqueness theorem]\label{thm:giut}
Let $\Lambda$ be a compactly aligned topological $k$-graph, and let $r^* : C_0(\Lambda^0)
\to C_0(\partial\Lambda)$ be the homomorphism $f \mapsto f\circ r$. Suppose that $\pi :
C^*(\Lambda) \to B$ is a homomorphism such that $\pi \circ r^*$ is injective on
$C_c(\Lambda^0)$. Suppose that there is a strongly continuous action $\beta : \TT^k \to
\Aut(B)$ such that for each $n \in \NN^k$ and $f \in C_c(\Gg_\Lambda)$ with $\supp(f)
\subseteq c^{-1}(n)$, we have $\beta_z(\pi(f)) = z^n \pi(f)$. Then $\pi$ is injective.
\end{theorem}

To prove the theorem, we first show that Yeend's boundary-path groupoid is amenable in the sense of
\cite{A-DR}, and then follow the standard argument of \cite{KP2000}. We begin by showing that the
kernel of $c$ is amenable. For us, amenability of $\Gg$ is important only as a hypothesis which
ensures that $C^*(\Gg_\Lambda)$ and $C^*_r(\Gg_\Lambda)$ coincide, so we will not dwell on the
rather technical definition. We thank Toke Carlsen for pointing out an error in an earlier version
of the proof of this result.

\begin{prop}\label{prp:kernel amenable}
Let $\Lambda$ be a compactly aligned topological $k$-graph. Then the kernel $H :=
c^{-1}(0)$ of $c$ is amenable, principal and satisfies $H^{(0)} = \Gg_\Lambda^{(0)}$.
\end{prop}
\begin{proof}
For $m \in \NN^k$, let $R_m$ denote the subgroupoid of $H$ defined by
\begin{align*}
R_m &:= \{(x,0,x) : x \in \partial\Lambda\} \cup \{(x,0,y) : d(x) = d(y) \ge m, \sigma^m(x) = \sigma^m(y)\} \\
    &= \{(x,0,x) : x \in \partial\Lambda\} \cup \{(\alpha z,0,\beta z) : z \in \partial\Lambda, \alpha, \beta \in \Lambda^m r(z)\}.
\end{align*}
Each $R_m$ is an equivalence relation, and each $R_m$ is also an $F_\sigma$ set (that is, a
countable union of closed sets) in $\partial\Lambda \times
\partial\Lambda$ because $\Gg_\Lambda$ is locally compact. We claim that each $R_m$ is proper as a
Borel groupoid \cite[Definition 2.1.2]{A-DR}. By \cite[Examples 2.1.4(2)]{A-DR}, this is equivalent
to the quotient space being a standard Borel space. The Mackey-Glimm-Ramsay dichotomy
\cite[Theorem~2.1]{Ramsay1990} implies that this in turn is equivalent here to the assertion that
the orbits are locally closed.

Fix $m \in \NN^k$. To see that the orbits in $R_m$ are indeed locally closed, first observe that
the orbit $[x]$ of $x$ in $R_m$ is equal to $\{x\}$ if $d(x) \not\ge m$, and is equal to $\{\alpha
\sigma^m(x) : \alpha \in \Lambda^m, s(\alpha) = x(m)\}$ otherwise. In the first case, $[x] = \{x\}$
is in fact closed because the topology on $\Gg_\Lambda^0$ is Hausdorff. In the second case, we
claim that
\begin{equation}\label{eq:orbit description}
[x] =
\Big(R_m^{(0)} \cap Z(\Lambda^m \ast_s \Lambda^m) \Big) \cap \overline{\{(\alpha \sigma^m(x), 0, \beta \sigma^m(x)) : d(\alpha) = d(\beta) = m\}}.
\end{equation}
To see this, observe that the right-hand side clearly contains $[x]$, so we need only show the
reverse inclusion. Fix
\[
(w,p,z) \in \Big(R_m^{(0)} \cap Z(\Lambda^m \ast_s \Lambda^m) \Big) \cap \overline{\{(\alpha \sigma^m(x), 0,
\beta \sigma^m(x)) : d(\alpha) = d(\beta) = m\}}.
\]
Then $w = z \in \partial\Lambda$, $p = 0$, and $d(w) \ge m$. Fix a sequence of pairs $(\alpha_j,
\beta_j) \in \Lambda^m x(m) \times \Lambda^m x(m)$ such that $(\alpha_j \sigma^m(x), 0,
\beta_j\sigma^m(x)) \to (w,0,w)$. In particular, $\alpha_j \sigma^m(x), \beta_j \sigma^m(x) \to w$
in $\partial\Lambda$. Then \cite[Proposition~3.12]{Y1}(i) ensures that the $\alpha_j
\sigma^m(x)(0,p \wedge d(x))$ converge to $w(0, m+p)$ for all $p \le d(w) - m$, which implies that
$d(w) \le d(x)$ and then that in fact $(\alpha_j \sigma^m(x))(0,q) \to w(0, q)$ for all $q \le
d(w)$. It therefore suffices to show that $d(w) \ge d(x)$. Suppose for contradiction that $i \le k$
satisfies $d(w)_i < d(x)_i$. Then $m_i \le d(w)_i < d(x)_i$, which gives
\[
(\alpha_j \sigma^m(x))(d(w), d(w) + e_i) = x(d(w), d(w) + e_i)\text{ for all $j \in \NN$.}
\]
In particular, the set
\[
J_{d(w), i} = \{j \in \NN : d(\alpha_j \sigma^m(x))_i \ge d(w) + e_i\}
\]
is equal to $\NN$ and hence infinite, but the sequence $(\alpha_j \sigma^m(x))(d(w), d(w) + e_i)_{j
\in J_{d(w), i}}$ is the constant sequence, and in particular is not wandering, contradicting
\cite[Proposition~3.12]{Y1}(ii). This proves~\eqref{eq:orbit description}.

Since $\Gg_\Lambda$ is \'etale, the unit space $R_m^{(0)} = \Gg_\Lambda^{(0)}$ is open. Each
$Z(\Lambda^n \ast_s \Lambda^n)$ is open in $\Gg_\Lambda$ by definition of the topology on
$\Gg_\Lambda$, so
\[
R_m^{(0)} \cap \bigcup_{n \ge m} Z(\Lambda^n \ast_s \Lambda^n)
\]
is open in the relative topology on $R_m$. Hence $[x]$ is the intersection of an open set and a
closed set in $R_m$ and hence is locally closed. So $R_n$ is a proper Borel groupoid; in particular
it is measurewise amenable.

The groupoid $H = \bigcup_{n \in \NN^k} R_n$ is therefore a direct limit (in the sense of
\cite[Section~5.3\textbf{f}]{A-DR}) of measurewise amenable groupoids, and hence is itself
measurewise amenable by \cite[Proposition~5.3.37]{A-DR}. Since $\Gg_\Lambda$ is \'etale by
\cite[Theorem~3.16 and Definition~4.8]{Y1}, $H$ is also \'etale. Since it is second-countable, it
follows that orbits are countable in $H$, so \cite[Theorem~3.3.7]{A-DR} implies that $H$ is
topologically amenable.

The factorisation property in $\Lambda$ implies that if $d(\alpha) = d(\beta)$, then for any $x \in
\partial \Lambda$, we have $\alpha x = \beta x$ if and only if $\alpha = \beta$. So $H$ is principal.
\end{proof}

In the earlier withdrawn version of this article, the first-, second- and fourth-named authors gave
an incorrect proof that if $\Gg$ is a second-countable, locally compact, Hausdorff, \'etale,
amenable groupoid and admits a continuous cocycle $c$ into an amenable group such that the kernel
of $c$ is an amenable groupoid, then $\Gg$ itself is amenable. Our proof was flawed because it
required strong surjectivity of $c$. The canonical $\ZZ^k$-valued cocycle on the groupoid of a
topological $k$-graph is usually not strongly surjective unless the range-map in $\Lambda$ is both
proper and surjective on each $\Lambda^n$, in which case Yeend's original results \cite{Y1} apply.
Fortunately, this gap in our argument is filled in by a recent result of Spielberg
\cite{Spielberg:xx11}.

\begin{cor}\label{thm:G_Lambda amenable}
Let $\Lambda$ be a compactly aligned topological $k$-graph, and let $\Gg_\Lambda$ be the
associated groupoid as in \cite{Y1}. Then $\Gg_\Lambda$ is (topologically) amenable and
$C^*(\Lambda)$ is nuclear.
\end{cor}
\begin{proof}
Proposition~9.3 of \cite{Spielberg:xx11} says that if $c$ is a continuous cocycle on a Hausdorff
\'etale groupoid $\Gg$ taking values in a discrete abelian group $G$ and the fixed-point groupoid
for $c$ is amenable, then $\Gg$ is amenable. To prove this result, Spielberg shows that $C^*(\Gg)$
is nuclear and then applies \cite[Corollary~6.2.14(ii)]{A-DR}. So the result follows from
Spielberg's argument combined with Proposition~\ref{prp:kernel amenable}.
\end{proof}

\begin{proof}[Proof of Theorem~\ref{thm:giut}]
Let $H := c^{-1}(0)$. Averaging over the gauge-action $\gamma : \TT^k \to \Aut(C^*(\Lambda))$
determines a faithful conditional expectation $\Phi^\gamma : C^*(\Gg_\Lambda) \to C^*(H)$.
Averaging over the action $\beta : \TT^k \to \Aut(B)$ determines a conditional expectation
$\Phi^\beta : B \to \pi(C^*(H))$ such that $\pi \circ \Phi^\gamma = \Phi^\beta \circ \pi$, so by a
standard argument it suffices to show that $\pi|_{C^*(H)}$ is injective.

By hypothesis $\pi$ is injective on $r^*(C_c(\Lambda^0))$, and it follows from
Proposition~\ref{prp:Lambda^0->G^0} that $\pi$ is injective on $C_0(\Gg_\Lambda^{(0)}) =
C_0(H^{(0)})$. Since Proposition~\ref{prp:kernel amenable} implies that $H$ is both principal and
amenable, it follows from~\cite[II, Proposition~4.6]{Renault1980} that $\pi|_{C^*(H)}$ is
injective.
\end{proof}

\section{The Cuntz-Krieger uniqueness theorem and simplicity}

In this section we use groupoid machinery to recover Yamashita's version of the Cuntz-Krieger
uniqueness theorem \cite{Yamashita:xx09}. (Yamashita's proof uses the technology of product systems
and Cuntz-Pimsner algebras.) We also use the results of \cite{BrownClarkEtAl:xx12} to characterise
simplicity of $C^*(\Lambda)$ in terms of the structure of $\Lambda$, and to establish a condition
under which $C^*(\Lambda)$ is purely infinite.

Recall from \cite{Y1} that given a topological higher-rank graph $\Lambda$, a boundary path $x \in
\partial \Lambda$ is said to be \emph{aperiodic} if $\sigma^m(x) \not= \sigma^n(x)$ for all
distinct $m,n \in \NN^k$ with $m,n \le d(x)$.

\begin{theorem}[The Cuntz-Krieger uniqueness theorem]\label{thm:CKUT}
Let $\Lambda$ be a compactly aligned topological $k$-graph. Let $r^* : C_0(\Lambda^0) \to
C_0(\partial\Lambda)$ be the homomorphism $f \mapsto f \circ r$. The following are
equivalent.
\begin{enumerate}
\item\label{it:ap} For every open set $V \subseteq \Lambda^0$ there exists an aperiodic element
    $x \in Z(V)$.
\item\label{it:CKUT} Every homomorphism $\pi : C^*(\Lambda) \to B$ such that $\pi \circ r^*$ is
    injective on $C_c(\Lambda^0)$ is an isomorphism.
\end{enumerate}
\end{theorem}

\begin{rmk}
Condition~(\ref{it:ap}) in Theorem~\ref{thm:CKUT} is precisely Yeend's aperiodicity condition~(A)
(see \cite[Theorem~5.2]{Y1}). Wright shows in Theorem~3.1 of \cite{Wright:xx11} that $\Lambda$
satisfies condition~(A) if and only if
\begin{equation}\label{eq:Sarah ap}
\parbox{0.9\textwidth}{for every pair $U,V$ of open subsets of $\Lambda$ such that $s(U) =
    s(V)$ and $s|_U, s|_V$ are homeomorphisms, there exists $\tau \in s(U)\Lambda$ such that
    $\MCE(U\tau, V\tau) = \emptyset$.}
\end{equation}
Since it does not involve elements of $\partial\Lambda$, which are hard to identify in practise,
this condition is easier to check in practice than condition~(\ref{it:ap}) of
Theorem~\ref{thm:CKUT} (see \cite[Section~4]{Wright:xx11}). We give an independent, although
somewhat round-about, proof of Wright's result in Lemma~\ref{lem:nlp vs BCFS} below.

The relationship between Yeend's aperiodicity condition and Yamashita's Condition~(B)
\cite[Definition~4.9]{Yamashita:xx09} is not transparent. However, since
\cite[Theorem~4.14]{Yamashita:xx09} says that Condition~(B) implies~(\ref{it:CKUT}) of
Theorem~\ref{thm:CKUT}, we deduce that Condition~(\ref{it:ap}) is at least formally weaker than
Condition~(B).
\end{rmk}

As in \cite{BrownClarkEtAl:xx12}, we say that a topological groupoid $\Gg$ is \emph{topologically
principal} if the set $\big\{u \in \Gg^{(0)} : \Gg^u_u = \{u\}\big\}$ of units with trivial
isotropy is dense in $\Gg^{(0)}$, and we say that $\Gg$ is \emph{minimal} if the only nonempty open
invariant subset of $\Gg^{(0)}$ is $\Gg^{(0)}$ itself.

\begin{proof}[Proof of Theorem~\ref{thm:CKUT}]
Lemma~\ref{lem:GLambda} says that $\Gg_\Lambda$ is second-countable, locally compact, Hausdorff and
\'etale. Corollary~\ref{thm:G_Lambda amenable} implies that it is amenable.

Theorem~5.2 of \cite{Y1} says that $\Lambda$ satisfies~(\ref{it:ap}) if and only if $\Gg_\Lambda$
is topologically principal.  Combined with the preceding paragraph,
\cite[Proposition~5.5]{BrownClarkEtAl:xx12} implies that $\Gg_\Lambda$ is topologically principal
if and only if every nontrivial ideal of $C^*(\Gg_\Lambda)$ has nontrivial intersection with
$C_0(\Gg_\Lambda^{(0)})$. Finally, Proposition~\ref{prp:Lambda^0->G^0} implies that every
nontrivial ideal of $C^*(\Gg_\Lambda)$ has nontrivial intersection with $C_0(\Gg_\Lambda^{(0)})$ if
and only if $C^*(\Lambda)$ satisfies~(\ref{it:CKUT}).
\end{proof}

We now employ the full strength of the characterisation \cite[Theorem~5.1]{BrownClarkEtAl:xx12} of
simplicity for $C^*$-algebras of second-countable locally-compact Hausdorff \'etale amenable
groupoids to characterise simplicity of topological higher-rank graph $C^*$-algebras.

\begin{theorem}\label{thm:simple}
Let $\Lambda$ be a compactly aligned topological $k$-graph. Then $C^*(\Lambda)$ is simple
if and only if both of the following conditions are satisfied:
\begin{enumerate}
\item $\Lambda$ satisfies condition~(\ref{it:ap}) of Theorem~\ref{thm:CKUT}; and
\item\label{it:cofinal} For every $x \in \partial\Lambda$ and open $U \subseteq \Lambda^0$
    there exists $n \in \NN^k$ such that $n \le d(x)$ and $U\Lambda x(n) \not= \emptyset$.
\end{enumerate}
\end{theorem}

\begin{lemma}\label{lem:nice open sets}
Let $\Lambda$ be a compactly aligned topological $k$-graph. Suppose that $V \subseteq
\Lambda^m$ is open, $F \subseteq \overline{V}\Lambda$ is compact and $Z(V) \cap Z(F)^c
\not= \emptyset$. Then there exists $p \ge m$ and a nonempty open subset $W$ of
$\Lambda^p$ such that $Z(W) \subseteq Z(V) \cap Z(F)^c$. In particular, if $U$ is an open
subset of $\partial\Lambda$, then there exist $n \in \NN^k$ and an open subset $W$ of
$\Lambda^n$ such that $Z(W) \subseteq U$.
\end{lemma}
\begin{proof}
We follow the argument of Theorem~5.2 of \cite{Y1}. Since $s|_{\Lambda^m}$ is a local
homeomorphism, we may assume that it restricts to a homeomorphism on $U$. Let $E := \{\lambda(m,
d(\lambda)) : \lambda \in F\}$. As in \cite[Definition~3.10]{Y1}, the map $\lambda \mapsto
\lambda(m,d(\lambda))$ is continuous on each $F \cap \Lambda^p$. Since $F$ is compact and $d :
\Lambda \to \NN^k$ is continuous, $d(F)$ is finite, and it follows that $E$ is compact. Fix $x \in
Z(V) \cap Z(F)^c$m and let $\lambda := x(0, m)$. Since $x \not\in Z(F)$, we have $\sigma^m(x)
\not\in Z(E)$. Since $\sigma^m(x) \in \partial\Lambda$ it follows that either $r(E)$ is not a
neighbourhood of $x(m)$, or $E$ is not exhaustive for $r(E)$. Suppose first that $r(E)$ is not a
neighbourhood of $x(m)$. Since $E$ is closed, it follows that there is an open neighbourhood $S$ of
$x(m)$ which does not intersect $r(E)$; and then $n := m$ and $W := US$ does the job. Now suppose
that $r(E)$ is not exhaustive for $r(E)$. Then there exists $\lambda \in r(E)\Lambda$ such that
$\lambda\Lambda \cap E\Lambda = \emptyset$. Since $E\Lambda$ closed there is then a neighbourhood
$S$ of $\lambda$ in $\Lambda^{d(\lambda)}$ such that $S\Lambda \cap E\Lambda = \emptyset$. Now $n
:= m + d(\lambda)$ and $W := US$ does the job.

The final statement follows because the $Z(V) \cap Z(F)^c$ are a base for the topology on
$\partial\Lambda$.
\end{proof}

\begin{lemma}\label{lem:cofinal vs irred}
Let $\Lambda$ be a compactly aligned topological $k$-graph. The following are equivalent:
\begin{enumerate}
\item $\Lambda$ satisfies condition~(\ref{it:cofinal}) of Theorem~\ref{thm:simple}.
\item $\Gg_\Lambda^{(0)}$ contains no nontrivial open invariant subsets.
\end{enumerate}
\end{lemma}
\begin{proof}
First suppose that $\Lambda$ satisfies condition~(\ref{it:cofinal}) of Theorem~\ref{thm:simple}.
Fix $x \in \partial\Lambda$. It suffices to show that $\overline{[x]} = \partial\lambda$. To see
this, fix $y \in \Gg_\Lambda^{(0)} = \partial\Lambda$. Each neighbourhood of $y$ contains a basic
open neighbourhood $Z(U) \cap Z(F)^c$ of $y$ where $U \subseteq \Lambda^m$ is relatively compact
and $F \subseteq \overline{U}\Lambda$ is compact. Lemma~\ref{lem:nice open sets} yields $p \in
\NN^k$ with $p \ge m$ and an open subset $W$ of $\Lambda^p$ such that $Z(W) \subseteq Z(U) \cap
Z(F)^c$. Proposition~4.3 of \cite{Y1} implies that each $v\partial\Lambda \not= \emptyset$, and so
$W\Lambda \cap F\Lambda = \emptyset$. Since $s(W)$ is open, condition~(\ref{it:cofinal}) of
Theorem~\ref{thm:simple} gives us $n \le d(x)$ such that $W\Lambda x(n) \not=\emptyset$, say
$\alpha \in W\Lambda x(n)$. Then $\alpha\sigma^n(x) \in \in [x] \cap Z(W) \subseteq [x] \cap Z(U)
\cap Z(F)^c$. Hence $y \in \overline{[x]}$.

Now suppose that $\Gg_\Lambda^{(0)}$ has no nontrivial open invariant subsets. Fix an open $U
\subseteq \Lambda^0$ and an element $x \in \partial\Lambda$. Since $\overline{[x]}$ is a nonempty
closed invariant set, it is all of $\partial\Lambda$. Since $Z(U)$ is open it follows that $[x]
\cap Z(U)$ is nonempty. By definition of $\Gg$, we have $[x] = \{\lambda\sigma^n(x) : n \in \NN^k,
\lambda \in \Lambda x(n)\}$, so $\Lambda$ satisfies condition~(\ref{it:cofinal}) of
Theorem~\ref{thm:simple}.
\end{proof}

\begin{proof}[Proof of Theorem~\ref{thm:simple}]
As in the proof of Theorem~\ref{thm:CKUT}, $\Gg_\Lambda$ is second-countable, locally compact,
Hausdorff, \'etale and amenable. Hence Theorem~5.1 of \cite{BrownClarkEtAl:xx12} implies that
$C^*(\Lambda)$ is simple if and only if $\Gg_\Lambda$ is topologically principal and minimal.
Theorem~5.1 of \cite{Y1} implies that $\Gg_\Lambda$ is topologically principal if and only if
$\Lambda$ satisfies condition~(\ref{it:ap}) of Theorem~\ref{thm:CKUT}. So the result follows from
Lemma~\ref{lem:cofinal vs irred}.
\end{proof}

Lemma~3.3 of \cite{BrownClarkEtAl:xx12} implies that a second-countable locally compact
Hausdorff groupoid $\Gg$ is topologically principal if and only if it satisfies the
apparently weaker condition (genuinely weaker in the absence of the assumption that $\Gg$
is second countable) that the interior of the isotropy subgroupoid $\bigcup_{u \in
\Gg^{(0)}} \Gg^u_u$ of $\Gg$ is precisely $\Gg^{(0)}$. Since this condition should be
easier to check, we describe what it says for a topological $k$-graph: it is a
topological analogue of the condition called ``no local periodicity" in \cite{RS07}. The
third condition below is Wright's finite-paths aperiodicity condition
\cite[Theorem~3.1(C)]{Wright:xx11}; as mentioned above, our argument below recovers the
equivalence (1)~$\iff$~(3) of \cite[Theorem~3.1]{Wright:xx11} via results of
\cite{BrownClarkEtAl:xx12} and \cite{Y1}.

Throughout the rest of the section we make frequent use of the notational convenience of
Notation~\ref{ntn:lazy}; that is, we write $Z(U *_s V) \cap Z(F)^c$ in place of $Z(U *_s V, p-q)
\cap Z(F,  p-q)^c$ when the former is unambiguous.

\begin{lemma}\label{lem:nlp vs BCFS}
Let $\Lambda$ be a compactly aligned topological $k$-graph. The following are equivalent:
\begin{enumerate}
\item\label{it:satisfies ap} $\Lambda$ satisfies condition~(\ref{it:ap}) of
    Theorem~\ref{thm:CKUT}.
\item\label{it:nlp} For every open set $V \subseteq \Lambda^0$ and every pair $m,n$ of distinct
    elements of $\NN^k$ there exists $x \in Z(V)$ such that either $d(x) \not\ge m \vee n$ or
    $\sigma^m(x) \not= \sigma^n(x)$.
\item\label{it:Sarah ap} $\Lambda$ satisfies~\eqref{eq:Sarah ap}.
\end{enumerate}
\end{lemma}
\begin{proof}
We first prove (\ref{it:satisfies ap})~$\iff$~(\ref{it:nlp}). We have seen that $\Lambda$ satisfies
condition~(\ref{it:ap}) of Theorem~\ref{thm:CKUT} if and only if $\Gg_\Lambda$ is topologically
principal. Lemmas 3.1~and~3.3 of \cite{BrownClarkEtAl:xx12} show that $\Gg_\Lambda$ is
topologically principal if and only if
\begin{equation}\label{eq:BCFS cond}
\parbox{0.9\textwidth}{every open subset of $\Gg_\Lambda
\setminus \Gg_\Lambda^{(0)}$ contains an element $(x,m,y)$ such that $x \not= y$.}
\end{equation}
So it suffices to show that~\eqref{eq:BCFS cond} is equivalent to~(\ref{it:nlp}).

First suppose that $\Lambda$ satisfies~(\ref{it:nlp}). Fix an open set $O \subseteq \Gg_\Lambda
\setminus \Gg_\Lambda^{(0)}$.

By definition of the topology on $\Gg_\Lambda$, the set $O$ contains a nonempty subset of the form
$Z(U *_s V) \cap Z(F)^c$ where $U \subseteq \Lambda^p$ and $V \subseteq \Lambda^q$ are relatively
compact with $s(U) = s(V)$, $s|_U$ and $s|_V$ are homeomorphisms, and $F$ is a compact subset of
$\bigcup_{\alpha \in \Lambda} \overline{U}\alpha \times \overline{V}\alpha$. Since the map
$(\mu\alpha,\nu\alpha) \mapsto \mu\alpha$ and the map $\mu\alpha \mapsto (\mu\alpha)(p, p +
d(\alpha))$ are continuous, it follows that $F = \{(\mu\alpha, \nu\alpha) : (\mu,\nu) \in U *_s V,
\alpha \in K\}$ for some compact $K \subseteq \overline{s(U)}\Lambda$. It suffices to find $(x,
p-q, y) \in Z(U *_s V) \cap Z(F)^c$ with $x \not= y$.

If $p = q$, then since $O \cap \Gg^{(0)}_\Lambda = \emptyset$ and since the groupoid $H$ of
Proposition~\ref{prp:kernel amenable} is principal, any $(x, 0, y) \in Z(U *_s V) \cap Z(F)^c$ does
the job. So we may suppose that $p \not= q$. By Lemma~\ref{lem:nice open sets}, there exists $m \in
\NN^k$ with $m \ge p$ and an open $W_0 \subseteq \Lambda^m$ such that $Z(W_0) \subseteq Z(U) \cap
Z(\overline{U}K)^c$. Let $n := m-p$ and let $W := \{\lambda(p,m) : \lambda \in W_0\} \subseteq
\Lambda^n$. Then $r(W) \subseteq s(U) = s(V)$ and $Z(UW *_s VW) \subseteq Z(U *_s V) \cap Z(F)^c
\subseteq O$. Let $p' := p - (p \wedge q)$ and $q' := q - (p \wedge q)$. Then $p \not= q$ forces
$p' \not= q'$. Since the source map in $\Lambda$ is open, $s(W) = s(W_0)$ is open, so
condition~(\ref{it:nlp}) implies that there exists $x \in Z(s(W))$ such that either $d(x) \not\ge
p' \vee q'$ or $\sigma^{p'}(x) \not= \sigma^{q'}(x)$. Let $\mu \in UW$ and $\nu \in VW$ be the
unique elements such that $s(\mu) = s(\nu) = r(x)$. Then $(\mu x, p-q, \nu x) \in O$.

We will show that $\mu x \not= \nu x$; equation~\eqref{eq:BCFS cond} then follows. We consider two
cases. First suppose that $d(x) \not\ge p' \vee q'$. Since $p' \wedge q' = 0$ it follows that there
exists $i \le k$ such that $d(x)_i < \infty$ and $p'_i \not= q'_i$. Thus $p_i \not= q_i$, and since
$d(\mu) = p + n$ and $d(\nu) = q + n$, it follows that $d(\mu)_i - d(\nu)_i \not= 0$. Since $d(x)_i
< \infty$, we have $d(\mu x)_i - d(\nu x)_i = d(\mu)_i - d(\nu)_i \not= 0$. In particular $d(\mu x)
\not= d(\nu x)$, forcing $\mu x \not= \nu x$ as required. Now suppose that $d(x) \ge p' \vee q'$.
Then~(\ref{it:nlp}) says that $\sigma^{p'}(x) \not= \sigma^{q'}(x)$. Since $\mu \in UW \subseteq
\Lambda^{p+n}$ and $\nu \in VW \subseteq \Lambda^{q+n}$, we have
\[
    \sigma^{p + n + q'}(\mu x) = \sigma^{q'}(x) \not= \sigma^{p'}(x) = \sigma^{q + n + p'}(\nu x).
\]
Since $p + n + q' = p + q - (p \wedge q) + n = q + n + p'$, we deduce that $\mu x \not= \nu x$ as
required.

Now suppose that $\Lambda$ does not satisfy~(\ref{it:nlp}). Fix an open set $V \subseteq \Lambda^0$
and distinct $m,n \in \NN^k$ such that $d(x) \ge m \vee n$ and $\sigma^m(x) = \sigma^n(x)$ for all
$x \in V\partial\Lambda$. Then $Z(V\Lambda^m *_s V\Lambda^n)$ is a nonempty open subset of
$\Gg_\Lambda$ which does not intersect $\Gg_\Lambda^{(0)}$ whose every element is an isotropy
element, and so~\eqref{eq:BCFS cond} does not hold. This completes the proof of (\ref{it:satisfies
ap})~$\iff$~(\ref{it:nlp})

We now establish (\ref{it:ap})~$\iff$~(\ref{it:Sarah ap}). As above, it suffices to show
that~\eqref{eq:BCFS cond} is equivalent to~(\ref{it:Sarah ap}).

First suppose that $\Lambda$ satisfies~\eqref{eq:Sarah ap}. Fix an open subset $B$ of $\Gg_\Lambda
\setminus \Gg_\Lambda^{(0)}$. As above there exist $m,n \in \NN^k$ and open sets $U \subseteq
\Lambda^m$ and $V \subseteq \Lambda^n$ such that $s(U) = s(V)$, $s|_U$ and $s|_V$ are
homeomorphisms and $Z(U *_s V) \subseteq B$. By~\eqref{eq:Sarah ap}, there exists $\tau \in
s(U)\Lambda$ such that $\MCE(U\tau, V\tau) = \emptyset$. Let $\alpha \in U$ and $\beta \in V$ be
the unique elements such that $s(\alpha) = s(\beta) = r(\tau)$ and fix $x \in s(\tau)
\partial\Lambda$. Then $g := (\alpha\tau x, m-n, \beta\tau x) \in Z(U *_s V) \subseteq B$, and
since $\MCE(\alpha\tau,\beta\tau) = \emptyset$, we have $\alpha\tau x \not= \beta\tau x$.

Now suppose that $\Lambda$ does not satisfy~\eqref{eq:Sarah ap}. So there exist $m,n \in \NN^k$ and
open $U \subseteq \Lambda^m$ and $V \subseteq \Lambda^n$ such that: (1) $s(U) = s(V) = W$, say; (2)
the source map restricts to homeomorphisms of $U$ and $V$ onto $W$; and (3) $\MCE(U\tau, V\tau)
\not= \emptyset$ for all $\tau \in W\Lambda$. By passing to subneighbourhoods, we may assume that
$\overline{U}$ and $\overline{V}$ are compact and contained in sets on which $s$ is a
homeomorphism, and that $\MCE(\overline{U}\tau, \overline{V}\tau) \not= \emptyset$ for all $\tau
\in s(\overline{U})$. Fix $\mu \in U$ and $\nu \in V$ with $s(\mu) = s(\nu)$. Then $\MCE(\mu,\nu)
\not= \emptyset$ (consider $\tau = s(\mu)$), so $\mu(0, m \wedge n) = \nu(0, m \wedge n)$ and for
each $\tau\in s(\mu)\Lambda$, we have
\[
    \mu(0, m\wedge n)\MCE(\mu(m\wedge n, m)\tau, \nu(m \wedge n, n)\tau) = \MCE(\mu\tau,\nu\tau) \not= \emptyset.
\]
So by replacing $U$ with $\{\mu(m \wedge n, m) : \mu \in U\}$ and $V$ with $\{\nu(m \wedge n, n) :
\nu \in V\}$, we may assume that $m \wedge n = 0$. We will show that $Z(U *_s V)$ consists entirely
of isotropy. We first establish the following claim.

\textbf{Claim.} \emph{For each $p \in \NN$, the set $\overline{W\Lambda^{pm}}$ is compact
exhaustive for each $v \in s(U)$.} The claim is trivial for $p = 0$, so suppose as an inductive
hypothesis that $\overline{W}\Lambda^{pm}$ is compact exhaustive for each $v \in s(U) = W$. Since
$m \wedge n = 0$ and hence $(p+1)m \wedge n = 0$, we have $\MCE(\overline{U\Lambda^{pm}},
\overline{V}) \subseteq \overline{V}\Lambda^{(p+1)m}$. Furthermore, for $\nu \in \overline{V}$ and
$\tau \in s(\mu) \Lambda^{(p+1)m}$, the element $\mu \in \overline{U}$ with $s(\mu) = r(\tau)$
satisfies $\MCE(\mu,\nu\tau) \not= \emptyset$, so $\nu\tau \in \MCE(\overline{U\Lambda^{pm}},
\overline{V})$. Hence
\[
    \MCE(\overline{U\Lambda^{pm}}, \overline{V}) = \overline{V}\Lambda^{(p+1)m}
    \quad\text{ for all $p \in \NN$.}
\]
Since each of $\overline{U}$ and $\overline{W}\Lambda^{pm}$ is compact, continuity of composition
implies that $\overline{U\Lambda^{pm}}$ is compact. Since $\overline{V}$ is compact also, and
$\Lambda$ is compactly aligned, it follows that $\overline{V}\Lambda^{(p+1)m}$ is compact. Since
$\lambda \mapsto \lambda(m, d(\lambda))$ is continuous on $\Lambda^{n+(p+1)m}$, we deduce that
$\overline{W}\Lambda^{(p+1)m}$ is compact. It remains to show that it is exhaustive for each $v \in
W$. For this fix $\tau \in W\Lambda$. The inductive hypothesis supplies an element $\eta$ of
$\MCE(\overline{W}\Lambda^{pm},\tau)$. By choice of $U$ and $V$, we have $\MCE(U\eta, V\eta) \not=
\emptyset$, say $\mu\eta\xi = \nu\eta\zeta \in \MCE(U\eta,V\eta)$ with $\mu \in U$ and $\nu \in V$.
By definition of $\eta$, we have $d(\eta) = (pm)\vee d(\tau) \ge pm$, so $d(\nu\eta\zeta) =
d(\mu\eta\xi) \ge d(\mu) + d(\eta) \ge (p+1)m$. Since $d(\nu) \wedge m = 0$, it follows that
$d(\eta\zeta) \ge (p+1)m$. Since $\eta(0, d(\tau)) = \tau$, we have $(\eta\zeta)(0, d(\tau) \vee
(p+1)m) \in \MCE(\tau, \overline{W}\Lambda^{(p+1)m})$. So $\overline{W}\Lambda^{(p+1)m}$ is
exhaustive for $v$. This proves the claim.

Now fix $(\mu,\nu) \in U *_s V$ and $x \in s(\mu)\partial\Lambda$ so that $(\mu x, m-n, \nu x)$ is
a typical element of $Z(U*_s V)$. We must show that $\mu x = \nu x$. The claim and the definition
of $\partial \Lambda$ imply that for each $p \in \NN$ there exists $\mu \in
\overline{W}\Lambda^{pm}$ such that $d(x) \ge d(\mu)$ and $x(0,pm) = \mu$. In particular, $d(x)_i =
\infty$ whenever $m_i > 0$, and similarly $d(x)_i = \infty$ whenever $n_i > 0$. So $d(\mu x) =
d(\nu x) = d(x)$, and $p \le d(x)$ if and only if $p \le d(\mu x)$. By choice of $U$ and $V$, we
have $\MCE(\mu x(0,p), \nu x(0,p)) \not= \emptyset$ for all $p \le d(x)$. Hence $(\mu x)(0, p) =
(\nu x)(0, p)$ for all $p \le d(\mu x) = d(\nu x)$. That is, $\mu x = \nu x$ as required.
\end{proof}

We use Anantharaman-Delaroche's criterion for pure infiniteness of a groupoid $C^*$-algebra
\cite[Proposition~2.4]{AD97} to provide a criterion under which $C^*(\Lambda)$ is simple and purely
infinite. Recall from \cite[Definition~2.1]{AD97} that a groupoid $\Gg$ is \emph{locally
contracting} if, for every open $U \subseteq \Gg^{(0)}$ there exist an open subset $V$ of $U$ and
an open bisection $B$ such that $\overline{V} \subseteq s(B)$ and $r(B\overline{V}) \subsetneq V$.

\begin{dfn}\label{dfn:contracting}
Given a compactly aligned topological $k$-graph $\Lambda$, we say that a precompact open
subset $U$ of $\Lambda^0$ is \emph{contracting} if there exist $m \not= n \in \NN^k$ and
nonempty precompact open sets $Y_m \subseteq \Lambda^m$ and $Y_n \subseteq \Lambda^n$
such that all of the following hold: $s(Y_m) = s(Y_n)$; $\overline{r(Y_m)} \subseteq
r(Y_n) = U$; the source map restricts to a homeomorphism on each of $Y_m$ and $Y_n$; for
every $\mu \in Y_m$ and $\nu \in Y_n$ such that $r(\mu) = r(\nu)$, we have $\MCE(\mu\tau,
\nu) \not= \emptyset$ for all $\tau \in s(\mu)\Lambda$; and there exists an open subset
$W$ of $Y_n\Lambda$ such that $\{\zeta(0,n) : \zeta \in W\} = Y_n$ and $\MCE(\mu, \zeta)
= \emptyset$ for all $\mu \in Y_m$ and $\zeta \in W$.
\end{dfn}

\begin{prop}\label{prp:locally contracting}
Let $\Lambda$ be a compactly aligned topological $k$-graph. Suppose that for every $v \in
\Lambda^0$ there exist $p \in \NN^k$ and an open set $V \subseteq \Lambda^p$ such that $v \in r(V)$
and $s(V)$ is contracting. Then $\Gg_\Lambda$ is locally contracting. If $\Lambda$ also satisfies
the hypotheses of Theorem~\ref{thm:simple}, then $C^*(\Lambda)$ is simple and purely infinite.
\end{prop}

To prove the proposition, we first prove that contracting neighbourhoods in $\Lambda^0$ give rise
to contracting bisections in $\Gg_\Lambda$.

\begin{lemma}\label{lem:technical}
Let $\Lambda$ be a compactly aligned topological $k$-graph. Suppose that $U \subseteq
\Lambda^0$ is contracting, and let $m$, $n$, $Y_m \subseteq \Lambda^m$, $Y_n \subseteq
\Lambda^n$ and $W \subseteq Y_n\Lambda$ be as in Definition~\ref{dfn:contracting}. Let
$Y_n'$ be a nonempty open set with $\overline{Y_n'} \subseteq Y_n$, and let $Y_m' :=
s^{-1}(Y'_n) \cap Y_m$. Then $\overline{r(Z(Y'_m *_s Y'_n))} \subsetneq s(Z(Y'_m *_s
Y'_n))$.
\end{lemma}
\begin{proof}
We first claim that $\MCE(Y_m, Y_n) = Y_m\Lambda^{(m \vee n) - m}$, and
$\overline{s(Y_m)}\Lambda^{(m \vee n) - m}$ is compact exhaustive for each $v \in s(Y_m)$.

The containment $\MCE(Y_m, Y_n) \subseteq Y_m\Lambda^{(m \vee n) - m}$ is clear. For the reverse
containment, fix $\tau \in s(Y_m)\Lambda^{(m \vee n) - m}$, let $\mu$ be the unique element of $Y_m
r(\tau)$, and fix $\nu \in Y_n$ such that $r(\nu) = r(\mu)$. By hypothesis, $\MCE(\mu\tau,\nu)
\not= \emptyset$, and since $d(\mu\tau) = m \vee n$, it follows that $\mu\tau \in \MCE(Y_m, Y_n)$.

To prove the claim, it remains to show that $\overline{s(Y_m)}\Lambda^{(m \vee n) - m}$
is compact exhaustive for each $v \in s(Y_m)$. First observe that $\MCE(\overline{Y_m},
\overline{Y_n})$ is compact because $\Lambda$ is compactly aligned. Since $\lambda
\mapsto \lambda(m, m \vee n)$ is continuous, it follows that $\{\lambda(m, m \vee n) :
\lambda \in \MCE(\overline{Y_m}, \overline{Y_n})\}$ is compact, so the first statement of
the lemma shows that $\overline{s(Y_m)}\Lambda^{(m \vee n) - m}$ is compact. To see that
it is exhaustive for each $v \in s(Y_m)$, fix $\tau \in s(Y_m)\Lambda$. Let $\mu \in Y_m$
and $\nu \in Y_n$ be the unique elements whose sources are equal to $r(\tau)$. By
hypothesis, we have $\MCE(\mu\tau, \nu) \not= \emptyset$, say $\mu\tau\alpha \in
\MCE(\mu\tau,\nu)$. Then $d(\mu\tau\alpha) \ge m \vee n$, and so $\eta := (\tau\alpha)(0,
(m \vee n) - m)$ belongs to $\overline{s(Y_m)}\Lambda^{(m \vee n) - m}$. In particular
$\tau \alpha \in \MCE(\tau, \eta) \subseteq \MCE(\tau, \overline{s(Y_m)}\Lambda^{(m \vee
n) - m})$. This proves the claim.

It follows from the claim and the definition of $\partial\Lambda$ that
\[
\overline{r(Z(Y'_m *_s Y'_n))}
    = \overline{Z(Y'_m)} \subseteq Z(Y'_n) =  s(Z(Y'_m *_s Y'_n)).
\]
To see that the containment is strict, observe that $Z(W) \cap Z(Y'_n)$ is a nonempty
open subset of $Z(Y'_n) \setminus Z(Y'_m)$.
\end{proof}

\begin{proof}[Proof of Proposition~\ref{prp:locally contracting}]
To see that $\Gg_\Lambda$ is locally contracting, fix a nonempty open subset $U$ of
$\Gg_\Lambda^{(0)} = \partial\Lambda$. By Lemma~\ref{lem:nice open sets} there exist $q
\in \NN^k$ and a nonempty open set $X \subseteq \Lambda^q$ such that $Z(\overline{X})
\subseteq U$. Since the source map is open, $s(X)$ is nonempty and open. Fix $v \in
s(X)$. By hypothesis, there exist $p \in \NN^k$ and an open set $V \subseteq \Lambda^p$
such that $v \in r(V)$ and $s(V)$ is contracting. Fix $m$, $n$, $Y_m$, $Y_n$ and $W$ as
in Definition~\ref{dfn:contracting}. Since each of $Y_m$, $Y_n$, $X$ and $V$ is open and
since composition is an open map, each of $X V Y_m$ and $XVY_n$ is open. Hence $Z(XVY_m)$
and $Z(XVY_n)$ are open. Let $B := Z(XVY_m *_s XVY_n)$. Then $B$ is a precompact open
bisection. Fix $\lambda \in X V Y_n$ and an open neighbourhood $Y_n'$ of $\lambda(p+q,
p+q+n)$ such that $\overline{Y'_n} \subseteq Y_n$. Let $K := Z(XVY'_n)$. Then
$\overline{K} \subseteq s(B)$, and Lemma~\ref{lem:technical} implies that
\begin{align*} r(B\overline{K})
    &= \overline{(\sigma^{p+q})^{-1}(r((Z(Y'_m *_s Y'_n)))) \cap Z(XV)} \\
    &= (\sigma^{p+q})^{-1}\big(\overline{r((Z(Y'_m *_s Y'_n)))}\big) \cap \overline{Z(XV)} \\
    &\subsetneq(\sigma^{p+q})^{-1}(s(Z(Y'_m *_s Y'_n))) \cap Z(XV)
    = K,
\end{align*}
so $\Gg_\Lambda$ is locally contracting as required.

For the final statement, observe that Theorem~5.1 of \cite{Y1} implies that $\Gg_\Lambda$ is
topologically principal and Lemma~\ref{lem:cofinal vs irred} implies that $\Gg_\Lambda$ is minimal.
Theorem~\ref{thm:G_Lambda amenable} implies that $\Gg_\Lambda$ is amenable. Hence
\cite[Proposition~2.4]{AD97} implies that $C^*(\Lambda) = C^*(\Gg_\Lambda)$ is simple and purely
infinite.
\end{proof}

\end{document}